\documentstyle[12pt]{article}
\newcommand{\sect}[1]{\section{#1}\setcounter{equation}{0}}

\font\mbn=msbm10 scaled \magstep1
\font\mbs=msbm7 scaled \magstep1
\font\mbss=msbm5 scaled \magstep1
\newfam\mbff
\textfont\mbff=\mbn
\scriptfont\mbff=\mbs
\scriptscriptfont\mbff=\mbss
\def\mbf{\fam\mbff}
\def\Re{{\mbf R}}

\def\Z{{\mbf Z}}
\def\Co{{\mbf C}}

\def\Di{{\mbf D}}
\newtheorem{Th}{Theorem}[section]

\newtheorem{C}[Th]{Corollary}
\newtheorem{D}[Th]{Definition}
\newtheorem{Proposition}[Th]{Proposition}
\newtheorem{R}[Th]{Remark}
\newtheorem{Ex}[Th]{Example}
\author{Alexander Brudnyi\thanks
{Research supported in part by NSERC.
}\\
Department of Mathematics\\
Ben-Gurion University of the Negev, Beer-Sheva\\
Israel}
\title{ON LOCAL BEHAVIOR OF ANALYTIC FUNCTIONS}
\date{}
\begin{document}
%==================================
%==================================
\maketitle
\begin{abstract}
We prove local inequalities for analytic functions defined on a convex
body in $\Re^{n}$ which generalize well-known classical inequalities 
for polynomials.
\end{abstract}
\sect{\hspace*{-1em}. Introduction.}
{\bf 1.1.} 
The classical Chebyshev inequality estimates the supremum norm of a 
univariate real polynomial $p$ on an interval $I$ by its norm on a subinterval
$I_{1}$ up to a multiplicative constant depending on the degree of $p$ and 
the ratio of lengths $|I|/|I_{1}|$ only.  
In the 1930's Remez [R] proved a generalization of the Chebyshev
inequality replacing $I_{1}$ by any measurable subset. A multivariate
inequality of such a kind (which coincides with the Remez inequality in the
one-dimensional case) was proved by Yu.Brudnyi and Ganzburg [BG] in
the 1970's. To formulate the result let 
${\cal P}_{k,n}(\Re)\subset\Re[x_{1},...,x_{n}]$ denote the space of real
polynomials of degree at most $k$ and $|U|$ denote the Lebesgue
measure of $U\subset\Re^{n}$. 
\\
{\bf Brudnyi-Ganzburg inequality}. {\em Let $V\subset\Re^{n}$ be a bounded
convex body and $\omega\subset V$ be a measurable subset. For every
$p\in {\cal P}_{k,n}$ the inequality}
\begin{equation}\label{bg}
\sup_{V}|p|\leq T_{k}\left(\frac{1+\beta_{n}(\lambda)}{1-\beta_{n}(\lambda)}
\right)\sup_{\omega}|p|
\end{equation}
{\em holds. Here
$\lambda:=|\omega |/|V|$ and $\beta_{n}=(1-\lambda)^{1/n}$ and
}
$T_{k}(x)=\frac{(x+\sqrt{x^{2}-1})^{k}+(x-\sqrt{x^{2}-1})^{k}}{2}$
{\em is the Chebyshev polynomial of degree $k$.} 
\\
From the above inequality one obtains
\begin{equation}\label{bg'}
\sup_{V}|p|\leq\left(\frac{4n|V|}{|\omega|}\right)^{k}\sup_{\omega}|p|\ .
\end{equation}
Inequalities of this kind are usually referred to as Bernstein-type
inequalities. They play an important role in the area of Approximation
Theory which investigates interrelation between analytic, approximative
and metric properties of functions. 
The purpose of this paper is to prove inequalities similar to
(\ref{bg'}) for analytic functions. We define  
the local degree of an analytic function which
expresses its geometric properties and generalizes the degree of a polynomial.
This notion is central for our consideration. It allows us  
to obtain better constants in Bernstein-type inequalities even in the 
standard (polynomial) case. 
It is worth pointing out that in recent years an essential progress was
done in studying Bernstein- and Markov-type inequalities for algebraic and
analytic functions. Such inequalities proved to be important in
different areas of modern analysis, see, e.g., [B], [BG], [BMLT], [Br1],
[Br2], [FN1], [FN2], [FN3], [G], [KY], [RY], [S]. We hope that the 
inequalities established in this paper would also have various applications 
in the fields that make use of the classical polynomial inequalities 
(Approximation Theory, trace and embedding theorems, signal processing, PDE 
etc). We proceed to formulation of the main result of the paper.

{\bf 1.2. A generalized Chebyshev inequality.}
Let $B_{c}(0,1)\subset B_{c}(0,r)\subset\Co^{n}$ be the pair of open complex
Euclidean balls of radii 1 and $r$ centered at $0$. Denote by
${\cal O}_{r}$ the set of holomorphic functions defined on $B_{c}(0,r)$.
Let $l_{x}\subset\Co^{n}(=\Re^{2n})$ be a real 
straight line passing through $x\in B_{c}(0,1)$. Further, 
let $I\subset l_{x}\cap B_{c}(0,1)$ be an interval and 
$\omega\subset I$ be a measurable subset.
\begin{Th}\label{re1}
For any $f\in {\cal O}_{r}$ there is a constant $d=d(f,r)>0$ such that
for any $\omega\subset I\subset l_{x}\cap B_{c}(0,1)$ 
\begin{equation}\label{rem1}
\sup_{I}|f|\leq\left(\frac{4|I|}{|\omega|}\right)^{d}\sup_{\omega}|f|\ .
\end{equation}
\end{Th}
\begin{Ex}\label{e1}
{\rm As an application of the above theorem we obtain local inequalities for
quasipolynomials.}
\begin{D}\label{quasi}
Let $f_{1},...,f_{k}\in (\Co^{n})^{*}$ be complex linear functionals.
A quasipolynomial with the spectrum $f_{1},...,f_{k}$ is a finite sum
$$
f(z)=\sum_{i}p_{i}(z)e^{f_{i}(z)},
$$
where $p_{i}\in\Co [z_{1},...,z_{n}]$. Expression $\sum_{i}(1+deg\ p_{i})$
is said to be the degree of $f$.
\end{D} 
\begin{Proposition}\label{qua}
Let $f$ be a quasipolynomial of degree $m$ and $l_{x}$ be a real straight
line passing through $x\in B_{c}(0,1)$. Then there is an absolute constant
$c>0$ such that the inequality 
$$
\sup_{I}|f|\leq\left(\frac{4|I|}{|\omega |}\right)^{c(\sqrt{k}M+m)}
$$
holds for any interval $I\subset l_{x}\cap B_{c}(0,1)$ and any measurable
subset $\omega\subset I$. Here $M:=\max_{i}\{||f_{i}||_{L^{2}(\Co^{n})}\}$.
\end{Proposition}
\end{Ex}
\begin{D}\label{degree}
The best constant $d$ in inequality (\ref{rem1}) will be called
the Chebyshev degree of the function $f\in {\cal O}_{r}$ in $B_{c}(0,1)$ 
and will be denoted by $d_{f}(r)$.
\end{D}
All constants in inequalities formulated below depend upon the possibilty to 
obtain an effective bound of Chebyshev degree in 
(\ref{rem1}). The following result gives such a bound 
in terms of the local geometry of $f$.

We say that a univariate holomorphic
function $f$ defined in a disk is ${\em p-valent}$ if it
assumes no value more than $p$-times there. We also say that $f$ is
$0$-valent if it is a constant. For any $t\in [1,r)$ let $L_{t}$ denote
the set of one-dimensional complex
affine spaces $l\subset\Co^{n}$ such that $l\cap B_{c}(0,t)\neq\emptyset$.
\begin{D}\label{d1}
Let $f\in {\cal O}_{r}$. The number
$$
v_{f}(t):=\sup_{l\in L_{t}}\{{\rm valency\ of}\ f|_{l\cap B_{c}(0,t)}\}
$$
is said to be the valency of $f$ in $B_{c}(0,t)$.
\end{D}
\begin{Proposition}\label{pr1}
For any $f\in {\cal O}_{r}$ and any $t$, $1\leq t<r$, the valency
$v_{f}(t)$ is finite. There is a constant $c=c(r)>0$ such that
$d_{f}(r)\leq cv_{f}(\frac{1+r}{2})$.
\end{Proposition}
\begin{R}\label{va}
For any holomorphic polynomial $p\in\Co[z_{1},...,z_{n}]$ of degree
at most $k$ the classical Remez inequality implies $d_{p}(r)\leq k$
while in many cases Proposition \ref{pr1} yields a sharper estimate.
\end{R}
 
{\bf 1.3.} In this section we formulate a generalization of inequality
(\ref{bg'}).
Let $B(0,1)\subset B_{c}(0,1)$ be the real Euclidean unit ball.
\begin{Th}\label{te4}
For any convex body 
$V\subset B(0,1)$, any measurable subset $\omega\subset V$
and any $f\in {\cal O}_{r}$ the inequality
$$
\sup_{V}|f|\leq\left(\frac{4n |V|}{|\omega|}\right)^{d_{f}(r)}
\sup_{\omega}|f|
$$
holds.
\end{Th}
The following corollary is a version of the $\log$-BMO-property for analytic
functions (cf. [St] and [Br2]).
\begin{C}\label{col1}
Under the hypothesis of Theorem \ref{te4} the inequality
$$
\frac{1}{|V|}\int_{V}\left|\log\frac{|f|}{||f||_{V}}\right|dx\leq Cd_{f}(r)
\log n
$$
holds with an absoulte constant $C$, where $||f||_{V}:=\sup_{V}|f|$.
\end{C}
Our next application of inequality (\ref{rem1}) is a 
generalization of Bourgain's polynomial inequality [B].
\begin{Th}\label{te5}
Let $V\subset B(0,1)$ be a convex body and 
$\widetilde d_{f}(r)$ be the smallest integer $\geq d_{f}(r)$. 
There are positive absolute constants
$c_{1},c_{2}$ such that the following inequality
\begin{equation}\label{b1}
|\{x\in V\ :\ |f(x)|>\frac{\lambda}{|V|}
\int_{V}|f(x)|dx\}|
\leq c_{1}\exp(-\lambda^{c_{2}/\widetilde d_{f}(r)})|V|
\end{equation}
holds for any $f\in {\cal O}_{r}$.
In particular,
$$
||f||_{L^{\Phi}(V,dx)}\leq (c_{1}+1)||f||_{L^{1}(V,dx)},\
$$
where $L^{\Phi}$ refers to the Orlicz space with the Orlicz function
$\Phi (t)=\exp(t^{c_{2}/\widetilde d_{f}(r)})-1$.
\end{Th}
\begin{R}\label{bou}
The original Bourgain's inequality for polynomials contains the degree of
the polynomial instead of $\widetilde d_{f}(r)$. 
\end{R}
As a corollary we also obtain the reverse H\"{o}lder inequality with the
constant which does not depend on the dimension (this result
does not follow from Theorem \ref{te4}).
\begin{C}\label{col2}
\begin{equation}\label{hold1}
\left(\frac{1}{|V|}\int_{V} |f(x)|^{s}dx\right)^{1/s}\leq
c(\widetilde d_{f}(r),s)
\frac{1}{|V|}\int_{V}|f(x)|dx\ \ \ \ \
(f\in {\cal O}_{r}, s\in\Z_{+})\ .
\end{equation}
\end{C}
The following example shows that in the polynomial case
our inequalities might be sharper than those of [BG] and [B].
\begin{Ex}\label{e2}
{\rm Let $f\in {\cal O}_{r}$ be such that $\sup_{B_{c}(0,r)}|f|<1$. Let
$\phi$ be a holomorphic non-polynomial function univalent in an open
neighbourhood $U$ of $\overline{\Di}=\{z\in\Co\ :\ |z|\leq 1\}$. 
Then using Proposition \ref{pr1} and Proposition \ref{propert} below yields
$d_{\phi\circ f}(r)\leq c(r)v_{f}(\frac{1+r}{2})$.
Consider a polynomial approximation $h_{k}$ of $\phi$ such that 
$deg\ h_{k}=k$ and $h_{k}$ is also univalent on $\Di$. 
Assume now that $f\in {\cal O}_{r}$ is a polynomial. Then
$deg(h_{k}\circ f)=k\cdot deg\ f$. Further, apply Brudnyi-Ganzburg
and Bourgain's polynomial inequalities to the polynomial $h_{k}\circ f$. Then 
the exponents in these inequalities will be equivalent to $k\cdot deg\ f$ and
$1/(k\cdot deg\ f)$, respectively.
However, in our generalizations of the above inequalities these
exponents contain numbers
$d_{h_{k}\circ f}(r)$ and $1/\widetilde d_{h_{k}\circ f}(r)$ with 
$d_{h_{k}\circ f}(r)\leq c(r)deg\ f$ and this is 
essentially better for all sufficiently large $k$.}
\end{Ex}
%=================================================
\sect{\hspace*{-1em}. Proofs of Theorem \ref{re1} and Proposition \ref{pr1}.}
{\bf 2.1.} We begin with auxiliary results used in the proof.\\
{\bf Parametrization of straight lines in the ball.} 
Let $B_{c}(0,s)$, $1<s<r$, be an open complex Euclidean ball. For any
$x\in B_{c}(0,s)$ consider the complex straight line
$l_{x,v}=\{x+vz\sqrt{s^{2}-|x|^{2}}\ ;\ \langle x,v\rangle =0,
|v|=1,\ z\in\Co\}$ passing through $x$. Here $|\cdot |$ denotes the 
Euclidean norm and $\langle .,.\rangle$ the inner product on $\Co^{n}$. 
In this way we 
parametrize the set $L_{s}$ of all complex straight lines passing through 
points of $B_{c}(0,s)$. Let
$f$ be a holomorphic function from ${\cal O}_{r}$. Consider the function
\begin{equation}\label{restrict}
F(z,x,v,s)=f(x+vz\sqrt{s^{2}-|x|^{2}})\ \ \ ( z\in\Di)
\end{equation}
Then $F(\cdot,x,v,s)$ is the restriction of $f$ to $l_{x,v}\cap B_{c}(0,s)$. 
Note also that for any $t<s$ the inequality
\begin{equation}\label{dilation}
\frac{s^{2}-|x|^{2}}{t^{2}-|x|^{2}}\geq (s/t)^{2}
\end{equation}
holds. This implies that the set $\{x+vz\sqrt{s^{2}-|x|^{2}}\ ;\ 
\langle x,v\rangle =0,\
|v|=1,\ z\in\frac{t}{s}\Di\}$ contains disk $l_{x,v}\cap B_{c}(0,t)$.
Set
$$
M(x,v,s,t)=\sup_{\frac{t}{s}\Di}\log |F(.,x,v,s)|\ .
$$
\begin{D}\label{index}
The number
$$
b_{f}(s,t,r):=\sup_{x,v}\{M(x,v,s,t)-M(x,v,s,1)\}
$$
is said to be the Bernstein index of $f\in{\cal O}_{r}$.
\end{D}
{\bf Bernstein index and Remez inequality.}
Assume that $F(\cdot,x,v,s)(=f|_{l_{x,v}\cap B_{c}(0,s)})$ has valency 
$m$ on $\frac{t}{s}\Di$. Assume also that $1<t<s$.
By Theorem 2.1.3 and Corollary 2.3.1 of [RY] (see
also [Br2, Lemma 3.1]), there is a constant $A=A(t)>0$ such that
\begin{equation}\label{royt}
M(x,v,s,(1+t)/2)-M(x,v,s,1)\leq Am\ .
\end{equation}
Then we apply the main inequality of Theorem 1.1
of [Br2] to the function $|F|$ obtaining that
there is a constant  $c=c(t,A)>0$ such that the inequality 
\begin{equation}\label{remroy}
\sup_{I'}|F|\leq\left(\frac{4|I'|}{|\omega'|}\right)^{cm}\sup_{\omega'}|F|
\end{equation}
is valid for any interval $I'\subset [-1/s,1/s]$ and any measurable 
set $\omega'\subset I'$.\\ 
Since $l_{x,v}\cap B_{c}(0,1)\subset
\{x+vz\sqrt{s^{2}-|x|^{2}}\ ;\ \langle x,v\rangle =0,\ |v|=1,\ z\in\frac{1}{s}\Di\}$,
(\ref{remroy}) implies inequality (\ref{rem1}) with exponent $cm$ for $f$ 
restricted to the real straight line $l_{x}\subset l_{x,v}$. 
%===========================================

{\bf 2.2. Proofs of Proposition \ref{pr1} and Theorem \ref{re1}.}
Let $1<t<r$ and $f\in {\cal O}_{r}$. First we prove inequality
$v_{f}(t)<\infty$. 

Fix a number $s$ satisfying $t<s<r$. For any $x\in B_{c}(0,s)$ consider 
complex straight
line $l_{x,v}=\{x+vz\sqrt{s^{2}-|x|^{2}}\ ;\ \langle x,v\rangle =0,\
|v|=1,\ z\in\Co\}$ passing through $x$.
Let $K:=\{(x,v)\in B_{c}(0,s)\times S^{2n-1}; \langle x,v\rangle =0\}$.
Further, for $f\in {\cal O}_{r}$ consider
the function $F$ defined by (\ref{restrict}).
Then $F$ is analytic on $\Di\times K$ and $F(\cdot,x,v,s)$ is holomorphic on 
$\Di$ for any $(x,v)\in K$. Let $K_{1}\subset K$ be a compact subset
that consists of points with the first coordinate from
$\overline{B_{c}(0,t)}$. In particular, the set of lines $l_{x,v}$ with
$x\in B_{c}(0,t)$ coincides with $L_{t}$ (defined just before 
Definition \ref{d1}).
Assume without loss of generality that $\sup_{B_{c}(0,s)}|f|=1$ and consider
the analytic function $\widetilde F(z,x,v,s,w)=F(z,x,v,s)-w$ defined
on $\Di\times K\times 2\Di$.  Set
$$
f_{1}(x,v,r,w)=\sup_{z\in\frac{2t}{t+s}\Di}\log |\widetilde F(z,x,v,s,w)|,\ \ \
f_{2}(x,v,r,w)=\sup_{z\in\frac{t}{s}\Di}\log |\widetilde F(z,x,v,s,w)|\ .
$$
Fix $(x,v,w)\in K_{1}\times\overline{\Di}$. If $\widetilde F(.,x,v,s,w)$
is not a constant then the number of its zeros in $\frac{t}{s}\overline{\Di}$
is estimated by the Jensen inequality
$$
\#\{z\in\frac{t}{s}\overline{\Di}\ :\ \widetilde F(z,x,v,s,w)=0\}\leq
c'(f_{1}(x,v,r,w)-f_{2}(x,v,r,w))
$$
with $c'=c'(s,t)>0$. Note also that by (\ref{dilation}),
the above number of zeros gives an upper bound for the number
of points $y\in l_{x,v}\cap B_{c}(0,t)$ such that $f(y)=w$.
Since $K_{1}\times\overline{\Di}$ is a compact, 
the Bernstein theorem of [FN3] and the 
Hadamard three circle theorem imply that
there is a constant $C=C(\widetilde F,K_{1}\times\overline{\Di})>0$ such that
$$
f_{1}(x,v,r,w)-f_{2}(x,v,r,w)\leq C
$$
for any $(x,v,w)\in K_{1}\times\overline{\Di}$. This inequality yields
$v_{f}(t)\leq c'C$ (see Definition \ref{d1}). 

It remains to prove
inequality $d_{f}(r)\leq c(r)v_{f}(\frac{1+r}{2})$. We will do it
in a parallel way with the proof of Theorem \ref{re1}.

Let $x\in B_{c}(0,1)$ and $l_{x}\subset\Co^{n}$ be a real straight line
passing through $x$. Let $I\subset l_{x}\cap B_{c}(0,1)$ be an
interval and $\omega\subset I$ be a measurable subset.
Set $s=\frac{1+r}{2}$, $t=\frac{1+s}{2}$ and
denote by $l_{x}^{c}=\{y+vz\sqrt{s^{2}-|y|^{2}}\ ;\ \langle y,v\rangle =0,\
|v|=1,\ z\in\Co\}$ the complex straight line containing $l_{x}$, where
$y\in l_{x}$ is such that $dist(0,l_{x})=|y|$.
By definition function
$F(.,y,v,s)=f|_{l_{x}^{c}\cap B_{c}(0,s)}$ determined by (\ref{restrict}) has 
valency $\leq v_{f}(s)$ on $\frac{t}{s}\Di$. Therefore Bernstein
index $b_{f}(s,\frac{1+t}{2},r)\leq Av_{f}(s)$ for $A=A(r)>0$ (see section
2.1). Finally, inequality (\ref{remroy}) and arguments of section 2.1
show that the inequality of Theorem \ref{re1} is valid with
$d\leq cv_{f}(s)$, $c=c(r)>0$. This implies that 
$$
d_{f}(r)\leq cv_{f}((1+r)/2)\ \ \ \ \ \Box  
$$ 
\begin{R}\label{bernst}
In order to estimate Chebyshev degree we can also use instead of 
$v_{f}(\frac{1+r}{2})$ an appropiriate Bernstein index 
$b_{f}(r)=b_{f}(s(r),t(r),r)$. Then
$d_{f}(r)\leq\widetilde{c}b_{f}(r)\leq cv_{f}(\frac{1+r}{2})$ with some
$\widetilde{c}=\widetilde{c}(r)>0$.
\end{R}
%=================================================
\sect{\hspace*{-1em}. Properties of Chebyshev Degree.}
We formulate further inequalities between Chebyshev degree and valency. 
In the following proposition the constant $c=c(r)$ is the same as in
Proposition \ref{pr1}.
\begin{Proposition}\label{propert}
(a)\ \ \ Let $f\in {\cal O}_{r}$ and $f(B_{c}(0,r))\subset\Di\subset\Co$.
Let $\phi$ be a holomorphic function defined in an open neighbourhood 
$U\supset\overline{\Di}$. Assume the $\phi$ has valency
$k$ in $U$. Then 
$$
d_{\phi\circ f}(r)\leq ckv_{f}((1+r)/2)\ .
$$
(b)\ \ \ Let $h:=e^{g}\in {\cal O}_{r}$. Then 
$$
d_{1/h}(r)\leq cv_{h}((1+r)/2)\ .
$$
(c)\ \ \ There is a constant $c_{1}=c_{1}(r)>0$ such that
$$
d_{fg}(r)\leq c_{1}(v_{f}((1+r)/2)+v_{g}((1+r)/2))
$$
for any $f,g\in {\cal O}_{r}$.
\end{Proposition}
Consider differential operator $(a,D)=\sum_{i=1}^{n}a_{i}D_{i}$,
where $a=(a_{1},...,a_{n})\in\Co^{n}$,
$D_{i}:=\frac{d}{dz_{i}}$, $i=1,...,n$ and
$z_{1},..., z_{n}$ are coordinates on $\Co^{n}$. 
Set $f_{m,a}:=(a,D)^{m}(f)$
\begin{Proposition}\label{Rol}
(The Rolle Theorem). Let $f\in {\cal O}_{r}$.
Assume that for any $a\in\Co^{n}$ the valency of $f_{m,a}$ satisfies
$v_{f_{m,a}}(\frac{1+3r}{4})\leq M$. Then there is a
constant $c_{2}=c_{2}(r)>0$ such that
$$
d_{f}(r)\leq c_{2}(m+M)\ .
$$
\end{Proposition}
%============================================
{\bf Proof of Proposition \ref{propert}.}
(a)\ According to the definition of the valency we have
$v_{\phi\circ f}(\frac{1+r}{2})\leq kv_{f}(\frac{1+r}{2})$, where $k$ is 
valency of $\phi$. Then $d_{\phi\circ f}(r)\leq
ckv_{f}(\frac{1+r}{2})$ by Proposition \ref{pr1}.\\
(b)\ The statement follows from Proposition \ref{pr1} and the identity
$v_{1/h}(\frac{1+r}{2})=v_{h}(\frac{1+r}{2})$ for $h=e^{g}$.\\
(c)\ According to results of section 2.1 it suffices to prove
the statement for univariate holomorphic functions 
$F(.,x,v,s)=f|_{l_{x,v}}$ and $G(.,x,v,s)=g|_{l_{x,v}}$. We consider more
general situation.

Assume that $\Di_{r_{1}}\subset\Di_{r_{2}}\subset\Co$, $r_{1}<r_{2}$, are
disks centered at 0 of radii $r_{1},r_{2}$, respectively. Further, assume
that $f,g$ are holomorphic in $\Di_{r_{2}}$ of valency $a$ and $b$,
respectively. We prove that
there is a constant $c=c(r_{1},r_{2})>0$ such that Chebyshev degree
$d_{fg}(r_{1})$ of $fg$ in $\Di_{r_{1}}$ $\leq c(a+b)$. Let
$K=\{z\in\Co\ :\ \frac{r_{1}+r_{2}}{2}\leq |z|\leq\frac{r_{1}+3r_{2}}{4}\}$
be an annulus in $\Di_{r_{2}}$ and
$$
g'=\frac{\log |g|-\sup_{\Di_{r_{2}}}\log |g|}{\sup_{\Di_{r_{2}}}\log |g|-
\sup_{\Di_{r_{1}}}\log |g|}\ .
$$
Repeating word-for-word the arguments of Lemma 2.3 of [Br2] we can find
a number $C=C(r_{1},r_{2})>0$ and a circle $S\subset K$ centered at $0$ such
that
$$
\inf_{S}g'\geq -C\
$$
(another relatively simple proof of this result can be done by Cartan's 
estimates for holomorphic functions, see, e.g. [L, p. 21]).
Going back to $|g|$ we obtain
$$
\inf_{S}|g|\geq\sup_{\Di_{r_{2}}}|g|\left(\frac{\sup_{\Di_{r_{1}}}|g|}{
\sup_{\Di_{r_{2}}}|g|}\right)^{C}\ .
$$
This implies
$$
\frac{\sup_{\Di_{r_{2}}}|fg|}{\sup_{S}|fg|}\leq\frac{\sup_{\Di_{r_{2}}}|f|
\sup_{\Di_{r_{2}}}|g|}{\sup_{S}|f|\inf_{S}|g|}\leq
\frac{\sup_{\Di_{r_{2}}}|f|}{\sup_{S}|f|}\cdot\left(
\frac{\sup_{\Di_{r_{2}}}|g|}{\sup_{\Di_{r_{1}}}|g|}\right)^{C}\ .
$$
Finally, according to Lemma 3.1 of [Br2] (see also section 2.1 above), there 
is a constant $B=B(r_{1},r_{2})>0$ such that
$$
\frac{\sup_{\Di_{r_{2}}}|f|}{\sup_{S}|f|}\leq
\frac{\sup_{\Di_{r_{2}}}|f|}{\sup_{\Di_{\frac{r_{1}+r_{2}}{2}}}|f|}
\leq B^{a}\ \ \ \ \
{\rm and}\ \ \ \ \
\frac{\sup_{\Di_{r_{2}}}|g|}{\sup_{\Di_{r_{1}}}|g|}\leq B^{b}\ .
$$
Thus we get
$$
\frac{\sup_{\Di_{r_{2}}}|fg|}{\sup_{\frac{r_{1}+3r_{2}}{4}}|fg|}\
\leq\frac{\sup_{\Di_{r_{2}}}|fg|}{\sup_{S}|fg|}\leq\widetilde{B}^{a+b},
$$
with $\widetilde{B}=\widetilde{B}(r_{1},r_{2}, B)>0$.
Then inequality (\ref{remroy}) applied to $|fg|$ implies the inequality of
Theorem \ref{re1} with exponent $c(a+b)$,
$c=c(r_{1},r_{2},\widetilde B)>0$. Therefore $d_{fg}(r_{1})\leq c(a+b)$.

In the multivariate case the above arguments estimate an appropriate
Bernstein index of $fg$ by sum of Bernstein indeces of $f$ and $g$.  
These indeces can be estimated by $c_{1}v_{f}(\frac{1+r}{2})$ and
$c_{1}v_{g}(\frac{1+r}{2})$ with some $c_{1}=c_{1}(r)>0$.
Thus according to Remark \ref{bernst},
$d_{fg}(r)\leq c'(r)(v_{f}(\frac{1+r}{2})+v_{g}(\frac{1+r}{2}))$. This
completes the proof of $(c)$. 

Proposition \ref{propert} is proved.\ \ \ \ \  $\Box$\\
{\bf Proof of Proposition \ref{Rol}.}
We, first, recall the relation between Bernstein index and Bernstein classes
(see [RY]).
\begin{D}\label{royo}
Let $f(z)=\sum_{i=0}^{\infty}a_{i}z^{i}$ be holmorphic in the disk $\Di_{R}$,
$R>1$.
We say that $f$ belongs to the Bernstein class $B_{N,R,c}^{2}$, if for
any $j>N$,
$$
|a_{j}|R^{j}\leq c\max_{0\leq i\leq N}|a_{i}|R^{i}\ .
$$
\end{D}
According to Corollary 2.3.1 of [RY], if the $m^{th}$ derivative $f^{(m)}$
of $f$ is $M$-valent then
$f^{(m+1)}\in B_{M-1,\frac{1+3R}{4},c^{M}}^{2}$ with $c:=c(R)>0$. Moreover,
from Definition \ref{royo} it follows that 
$f\in B_{m+M,\frac{1+3R}{4},c^{M}}^{2}$.
Then Theorem 2.1.3 of [RY] based on the last implication yields
\begin{equation}\label{plus}
\sup_{\Di_{\frac{1+R}{2}}}|f|\leq
a^{m+M}{\sup_{\Di_{1}}|f|}
\end{equation}
for some constant $a=a(R)>1$.

We proceed with the proof of the proposition. As in the proof of 
Proposition \ref{propert} it suffices to prove the result for
restriction $F_{l}$ of $f$ to a complex line $l$ passing through a point 
of $B_{c}(0,1)$. Then the condition of the proposition implies that
$m^{th}$ derivative $F_{l}^{(m)}$ of $F_{l}$ has valency 
at most $M$ in the larger disk $l\cap B_{c}(0,\frac{1+3r}{4})$.
Therefore the required result follows immeadiately from inequality 
(\ref{plus}) 
(an estimate for Bernstein index) and arguments of section 2.1. 

The proof of proposition is complete.\ \ \ \ \ $\Box$\\
%====================================
\sect{\hspace*{-1em}. Proofs.}
\noindent {\bf Proof of Proposition \ref{qua}.}
Let $l_{y}^{c}=
\{y+vz\sqrt{4-|y|^{2}}\ ;\ \langle y,v\rangle =0,\ |v|=1,\ z\in\Co\}$
be a complex straight line passing through a point $y\in B_{c}(0,1)$.
Consider restriction $F$ of the quasipolynomial
$f(z)=\sum_{i=1}^{k}p_{i}(z)e^{f_{i}(z)}$ to $l_{y}^{c}$. Then $F$ is
a univariate quasipolynomial of the form
$$
F(z)=q_{i}(z)e^{f_{i}(y)}e^{z\sqrt{4-|y|^{2}}f_{i}(v)}\ \ \ \ \
(q_{i}\in\Co[z])
$$
of degree $\leq m$. 
We estimate valency of $F$ in disk $\Di_{2}:=2\Di$ (i.e.
we estimate the number of zeros of $F+c$ for any $c\in\Co$). 
Note that $F+c$ is also a
quasipolynomial of degree $\leq m+1$. Further, by definition
$\max_{i}\{|f_{i}(v)|\}\leq M$ implying
$\sqrt{4-|y|^{2}}f_{i}(v)\in\Di_{2M}$ for any $i$. Then by
Theorem 2 in [KY] the number of zeros of
$F+c$ in $\Di_{2}$ less than or equal to 
$m+\frac{2}{\pi}(\sqrt{k+1}+1)\cdot 16M<32(\sqrt{k+1}M+m)$. 
This and Proposition \ref{pr1} yield 
$$
d_{f}(2)\leq cv_{f}(3/2)\leq c'(\sqrt{k+1}M+m)
$$ 
with an absolute constant $c'>0$. The required inequality follows from
the definition of Chebyshev degree.
\ \ \ \ \ $\Box$\\
{\bf Proof of Theorem \ref{te4}.}
Let $V\subset B(0,1)$ be a convex body, $\lambda\subset V$ be a measurable
subset and $f\in {\cal O}_{r}$. Take a point $x\in V$ such that
$$
|f(x)|=\sup_{V}|f|\ .
$$
(Without loss of generality we may assume that $x$ is an interior point of
$V$; for otherwise, apply the arguments below to an interior point 
$x_{\epsilon}\in V$, $\epsilon>0$, such that 
$|f(x_{\epsilon})|>\sup_{V}|f|-\epsilon $ and then take the limit when 
$\epsilon\to 0$.)
According to Lemma 3 of [BG] there is a ray $l$ with origin at 
$x$ such that
\begin{equation}\label{geom}
\frac{mes_{1}(l\cap V)}{mes_{1}(l\cap\lambda)}\leq\frac{n|V|}{|\lambda|}.
\end{equation}
Let $l'$ be the real straight line containing $l$. Applying inequality
(\ref{rem1}) to $f|_{l'}$ with $I:=l\cap V$ and
$\omega:=l\cap\lambda$ and then inequality (\ref{geom}) lead to the 
required result.\ \ \ \ \ $\Box$
\begin{R}\label{doubl}
Assume that $\omega\subset V$ is a pair of Euclidean balls of radii $R_{1}$ 
and $R_{2}$, respectively. Then the ray $l$ in (\ref{geom}) can be chosen 
such that the constant in the inequality of Theorem \ref{te4} will be  
$\left(\frac{4R_{1}}{R_{2}}\right)^{d_{f}}$.
\end{R}
{\bf Proof of Corollary \ref{col1}.}
Let $V\subset B(0,1)$ be a convex body and $f\in {\cal O}_{r}$. For the
distribution function $D_{f}(t):=mes\{x\in V : |f(x)|\leq t\}$ the inequality
of Theorem \ref{te4} acquires the form
$$
D_{f}(t)\leq 4n |V|\left(\frac{t}{||f||_{V}}\right)^{1/d_{f}(r)}.
$$
The required result follows from the above inequality and the identity
$$
\int_{V}\left|\log\frac{|f|}{||f||_{V}}\right|dx
=\int_{0}^{|V|}\left|\log\frac{f_{*}}{||f||_{V}}\right|dx\ ,
$$
where $f_{*}=\inf\{s\ :\ D_{f}(s)\geq t\}$ (cf. [Br1, Th. 5.1]).
\ \ \ \ \  $\Box$\\
{\bf Proof of Theorem \ref{te5}.}
Let $V\subset B(0,1)$ be a convex body. For a real straight line
$l$, $l\cap V\neq\emptyset$, and an interval $I\subset l\cap V$
inequality (\ref{rem1}) implies 
$$
mes\{t\in I\ :\ |f(t)|\geq 10^{-\widetilde d_{f}(r)}||f||_{I}\}>|I|/2
$$
holds for any $f\in {\cal O}_{r}$ with $||f||_{I}=\sup_{I}|f|$.
Applying the same arguments as in the original
proof of Bourgain's inequality for polynomials [B] but based on the above 
inequality instead of that of Lemma 3.1 of [B] one obtains the required 
result.
The second part of Theorem \ref{te5} follows from the distributional 
inequality of the theorem and the definition
$$
||f||_{L^{\Phi}(V,dx)}:=\inf\{A\geq 0\ :\ \int_{V}\Phi(|f|/A)dx\leq 1\}
\ \ \ \ \ \Box
$$
{\bf Proof of Corollary \ref{col2}.}
The reverse H\"{o}lder inequality (\ref{hold1}) follows 
straightforwardly from the
distributional inequality of Theorem \ref{te5}.\ \ \ \ \ $\Box$
%===========================================
\sect{\hspace*{-1em}. Concluding Remarks.}
{\bf 5.1.}
Consider a uniformly bounded sequence of functions $\{f_{i}\}_{i\in I}
\subset {\cal O}_{r}$ and define $h_{i}=(|f_{i}|)^{1/d_{f_{i}}(r)}$. Let
$$
h=(\overline{\lim_{i\to\infty}}h_{i})^{*},
$$
where $g^{*}$ denotes upper semicontinuous regularization of $g$. Clearly
$h$ is logarithmically plurisubharmonic.
Then one can show that
inequalities of Theorems \ref{te4} and \ref{te5} hold for $h$ with
exponents $1$ and $c_{2}$ instead of $d_{f}(r)$ and 
$c_{2}/\widetilde d_{f}(r)$,
respectively.

Assume that a plurisubharmonic function $u$ is taken from the class $L$, 
i.e. satisfies
$$
u(z)\leq\alpha +\log(1+|z|)\ \ \ \ \  (z\in\Co^{n})
$$
for some $\alpha\in\Re$. Then inequalities of Theorems \ref{te4} and 
\ref{te5} are valid for $e^{u}$ restricted to a convex body $V\subset\Re^{n}$
with the constants which contain exponents 1 and  $c$ (absolute constant),
respectively. It follows from the fact $\displaystyle
u=(\overline{\lim_{i\to\infty}}(\log|p_{i}|)/deg\ p_{i})^{*}$,
where $\{p_{i}\}$ is a sequence of holomorphic polynomials on $\Co^{n}$
(for the proof see, e.g. [K]).

{\bf 5.2.}
Inequalities of Theorems \ref{te4} and \ref{te5} can also be written in the
same form for convex bodies in $B_{c}(0,1)$, where one 
replaces coefficient $4n$ by $8n$ in the first inequality.

{\bf 5.3.} If $f_{1},...,f_{k}$ are functions from ${\cal O}_{r}$ and
$p$ is a holomorphic polynomial of degree $d$ then for $h=p(f_{1},...,f_{k})$ 
its degree
$d_{h}(r)$ is bounded by a constant depending on $d,r$ and $f_{1},...,f_{k}$.
It follows, e.g., from results of [FN3] and arguments used in the proof of
Proposition \ref{pr1}. However, it is difficult to
obtain an explicit estimate for $d_{h}(r)$ even in the case of naturally
defined functions $f_{i}$ (e.g., taken as solutions of some systems of ODEs).
Assume, e.g., that $f_{1}=z_{1},...,f_{n}=z_{n}$ are coordinate functions on 
$\Co^{n}$ and $k\geq n$. Then inequality $d_{h}(r)\leq c d$ holds
for any polynomial $p$ of degree $d$ with $c$ which does not depend on $d$ if 
and only if $f_{n+1},...,f_{k}$
are algebraic functions, see [S] and [Br2, Th.1.3].

\end{document}